\input amstex
\documentstyle{amsppt}
\magnification=\magstep1                        
\hsize6.5truein\vsize8.9truein                  
\NoRunningHeads
\loadeusm

\magnification=\magstep1                        
\hsize6.5truein\vsize8.9truein                  
\NoRunningHeads
\loadeusm

\document
\topmatter

\title
Bourgain-Ruzsa type estimates for the maximum and minimum of Littlewood Cosine polynomials on the real line
\endtitle

\rightheadtext{ultraflat unimodular polynomials}

\author Tam\'as Erd\'elyi
\endauthor

\address Department of Mathematics, Texas A\&M University,
College Station, Texas 77843, College Station, Texas 77843 \endaddress

\endthanks


\date May 12, 2026
\enddate

\email terdelyi\@tamu.edu
\endemail


\abstract
We prove that 
$$\max_{t \in {\Bbb R}}{T_n(t)} \geq \frac{1}{60}n^{1/3} \qquad \text {\rm  and} \qquad -\min_{t \in {\Bbb R}}{T_n(t)} \geq \frac{1}{60}n^{1/3}$$
for every trigonometric polynomial $T_n$ of the form 
$$T_n(t) = \sum_{j=1}^n{\cos(jt-\theta_j)}\,, \quad \theta_j \in {\Bbb R}\,, \quad t \in {\Bbb R}\,.$$
\endabstract
\endtopmatter

\head 1. Introduction \endhead
Let $0 = \lambda_0 < \lambda_1 < \cdots < \lambda_m$. Our starting point is the following question.
How large can the maximum of a trigonometric polynomial
$$S_m(t) = \sum_{j=0}^m{A_j \cos(\lambda_jt)}\,, \qquad A_j \in {\Bbb R}\,,$$
be on the real line? Since
$$\int_{0}^{2\pi}{|S_m(t)|^2\, dt} = 2\pi \left( |A_0|^2 + \frac 12 \sum_{j=1}^m{|A_j|^2} \right)\,,$$
the inequality
$$\max_{t \in [0,2\pi]}{|S_m(t)|} \geq \left( |A_0|^2 + \frac 12 \sum_{j=1}^m{|A_j|^2} \right)^{1/2}$$
obviously holds. But how large can
$$\max_{t \in [0,2\pi]}{S_m(t)} \tag 1.1$$
be? To give a decent lower bound for (1.1) looks rather difficult.

The result below stated in [BE07] is straightforward from [DL93, pages 285-288] which offers an
elegant book proof of the Littlewood Conjecture first shown in [Ko81] and [Mc81]. The book
[Bo02] deals with a number of related topics. Littlewood [Li61, Li64, Li66, Li68] was interested in many closely
related problems.
	
\proclaim{Theorem 1.1} Let $\lambda_0 < \lambda_1 < \cdots < \lambda_m$ be nonnegative integers and let
$$S_m(t) = \sum_{j=0}^m{A_j \cos(\lambda_jt)}\,, \qquad A_j \in {\Bbb R}\,.$$
Then
$$\int_{0}^{2\pi}{|S_m(t)| \, dt} \geq \frac{1}{60} \sum_{j=0}^m{\frac{|A_{m-j}|}{j+1}}\,.$$
\endproclaim

The above theorem can be used to obtain a nontrivial lower bound for (1.1).

\proclaim{Theorem 1.2} Let $\lambda_1 < \lambda_2 < \cdots < \lambda_m$
be positive integers and let
$$S_m(t) = \sum_{j=1}^m{A_j \cos(\lambda_jt)}\,, \qquad A_j \in {\Bbb R}\,.$$
Then
$$\max_{t \in [0,2\pi]}{S_m(t)} \geq \frac{1}{120} \sum_{j=0}^{m-1}{\frac{|A_{m-j}|}{j+1}}\,.$$
\endproclaim

See also the Appendix of [Bo93] as a related topic. To see Theorem 1.2 observe that
$$\int_{0}^{2\pi}{S_m(t) \, dt} = 0\,,$$
and hence with $S_m^+(t) := \max\{S_m(t),0\}$ and $S_m^-(t) := \min\{S_m(t),0\}$, we have
$$\int_{0}^{2\pi}{S_m^+(t) \, dt} = \int_{0}^{2\pi}{S_m^-(t) \, dt} = \frac 12 \int_{0}^{2\pi}{|S_m(t)| \, dt}\,,$$
which, together with Theorem 1.1, gives Theorem 1.2.

In [Er11] Theorem 1.2 is significantly improved for the interesting classes of Littlewood cosine polynomials
$$T_q(t) = \sum_{j=0}^q{a_j\cos(jt)}, \qquad  a_j \in \{-1,1\}\,,$$
at least in the case when $p = 2q+1$ is an odd prime where we rely heavily on Ruzsa's 
paper [Ru04], who claims the best result today to solve Chowla's Cosine Problem in [Ch65] below.

\proclaim {Chowla's Cosine Problem}
Let $A \subset {\Bbb N}$ be a finite set of distinct integers and set
$$\mu(A) := -\min_{t \in [0,2\pi]}{\sum_{a \in A}{\cos(at)}}\,.$$
What is $\nu(n) := \min{\{\mu(A): A \subset {\Bbb N}, \, |A| = n\}}$?
\endproclaim

In the Introduction of [Ru04] Ruzsa writes:
``Let $A$ be a finite set of positive integers, $|A| = n$, and write
$$f(x) = \sum_{a \in A}{\cos(ax)}\,.$$ 
Since $f(0) > 0$ and $\int_0^{2\pi}{f(x)\,dx} = 0$,
we have $\min{f(x)} < 0$.
It is a difficult question to estimate this minimum uniformly for every set of size $n$.
Bourgain [Bo84] proved
$$\min{f(x)} < -c_1 \exp(c_2(\log n)^{c_3})$$
with unspecified absolute constants $c_1$, $c_2$, and $c_3$. In another paper [Bo86]
he showed that one can take $c_3 = 1/2$ under the assumption that $A \subset [1,n2^{\sqrt{\log n}}]$.
Our aim is to prove this without restriction.

\proclaim{Theorem A}
With the above notations we have
$$\min{f(x)} < -c_4 \exp(c_5(\log n)^{1/2})$$
with a positive absolute constant $c_4$ and $c_5 = \sqrt{(\log 2)/8}$.
\endproclaim

\noindent $\ldots$". 

Note that $\min{f(x)}$ in the above quotation denotes the smallest value
of $f(x)$ on the real number line ${\Bbb R}$. Note also that the above quotation
corrects two misprints in Ruzsa's paper by removing the minus sign from the exponent
at two places.

We denote the additive group of $p$ elements $\{0,1,\ldots,p-1\}$ under addition modulo $p$
by ${\Bbb Z}_{p}$. Let $y_j := j/p$ for $j = 0,1, \ldots, p-1$,
$$E_p := \{y_0, y_1, \ldots, y_{p-1}\} \qquad \text {\rm and} \qquad
E_p^* := E_p \cup \textstyle {\left\{ \frac{3}{2p} \right\}}\,.$$
In [Er11] the following result is proved.

\proclaim{Theorem 1.3} If $p = 2q + 1$ is a prime, then the maximum of a $T_q$ defined by
$$T_q(2\pi t) := \sum_{j=0}^q{a_j\cos(2\pi jt)}\,, \qquad  a_j \in \{-1,1\},$$
on $E_p^*$ is at least $c_1 \exp(c_2(\log q)^{1/2})$,
with an absolute constant $c_1 > 0$ and $c_2 = \sqrt{(\log 2)/8}$.
\endproclaim

\head 2. New Results \endhead

Let
$${\Cal LC}_n := \left\{ T_n: T_n(t) = \sum_{j=1}^n{\cos(jt + \theta_j)}\,, \quad \theta_j \in {\Bbb R}\,, \quad t \in {\Bbb R} \right\}\,.$$
Note that the elements of the class 
$${\Cal LC}^*_n := \left\{ T_n: T_n(t) = \sum_{j=1}^n{a_j\cos(jt)}\,, \quad a_j \in \{-1,1\}\,, \quad t \in {\Bbb R} \right\}$$
are in ${\Cal LC}_n$, as it can be see by choosing $\theta_j = 0$ if $a_j = 1$ and by choosing $\theta_j = \pi$ if $a_j = -1$.

\proclaim{Theorem 2.1}
We have  
$$\max_{t \in {\Bbb R}}{T_n(t)} \geq \frac{1}{60}n^{1/3}$$
for every polynomial $T_n \in {\Cal LC}_n$.
\endproclaim

\proclaim{Corollary 2.2} 
We have  
$$-\min_{t \in {\Bbb R}}{T_n(t)} \geq \frac{1}{60}n^{1/3}$$
for every polynomial $T_n \in {\Cal LC}_n$.
\endproclaim

\proclaim{Remark 2.3}
We have 
$$\max \left( \max_{t \in {\Bbb R}}{T_n(t)},  -\min_{t \in {\Bbb R}}{T_n(t)} \right) \geq \left(n/4\right)^{1/2}$$
for every $T_n \in {\Cal LC}_n$.
\endproclaim

\head 3. Lemmas \endhead

To prove Theorem 2.1 we need the Bernstein inequality in $L_q$ for trigonometric polynomials in $L_q$, $q > 0$.
In fact, we need only the case $q=1$ of it. See [Ar79] and [Ar82] for the general case $q > 0$. For a book proof see 
[DL93] or [Er20]. For the case $q \geq 1$ see [BE95]. 

\proclaim{Lemma 3.1} We have 
$$\int_{0}^{2\pi}{\left|T_n^{\prime}(t)\right|^q \, dt} \leq n^q \int_{0}^{2\pi}{\left|T_n(t)\right|^q \, dt}$$
for every trigonometric polynomial $T_n$ of the form 
$$T_n(t) = \sum_{j=-n}^n{a_je^{ijt}}\,, \qquad a_j \in {\Bbb C}\,.$$
\endproclaim

We call the functions 
$$T_n(t) = a_0 + \sum_{k=1}^n{a_k \cos(kt) + b_k \sin(kt)}\,, \qquad a_k, b_k \in {\Bbb R}, \quad a_nb_n \neq 0\,,$$
a real trigonometric polynomial of degree $n$. The following inequality is a version of the Bernstein-Szeg\H o inequality.
For a proof of it see [BE95]. We will use it in the proof of Theorem 2.1 as well.

\proclaim{Lemma 3.2}
Suppose $T_n$ is a real trigonometric polynomial of degree at most $n$ satisfying $m \leq T_n(t) \leq M$ for all $t \in {\Bbb R}$.
We have 
$$\left|T_n^{\prime}(t)\right|^2 + n^2\left|T_n(t) - \frac{M+m}{2}\right|^2 \leq n^2\left(\frac{M-m}{2}\right)^2\,, \qquad t \in {\Bbb R}\,.$$
\endproclaim

\head 4. Proofs \endhead

\demo{Proof of Theorem 2.1}
Let $T_n \in {\Cal LC}_n$. Let $T_n$ be a trigonometric polynomial of the form
$$T_n(t) = \sum_{j=1}^n{\cos(jt + \theta_j)}\,, \quad \theta_j \in {\Bbb R}\,, \quad t \in {\Bbb R}\,. \tag 4.1$$
Therefore by Parseval's formula we get that
$$\int_{0}^{2\pi}{\left|T_n^{\prime}(t)\right|^2 \, dt} = \pi \sum_{j=1}^{n}{j^2} = \pi \frac{n(n+1)(2n+1)}{6} \tag 4.2$$  
We have
$$\int_E{\left|T_n^{\prime}(t)\right|^2 \, dt} \leq \left( \int_{0}^{2\pi}{|T_n^{\prime}(t)| \, dt} \right) \left( \max_{t \in E}{|T_n^{\prime}(t)|} \right) \tag 4.3$$ 
for every measurable set $E \subset {\Bbb R}$.
Suppose that there are $0 \delta_n \leq n/16$ such that 
$$\max_{t \in {\Bbb R}}{|T_n(t)|} \leq \delta_n\,. \tag 4.4$$ 
As
$$\int_{0}^{2\pi}{T_n(t) \, dt} = 0$$
we have
$$\int_{0}^{2\pi}{T_n^+(t) \, dt} = \int_{0}^{2\pi}{T_n^-(t) \, dt}\,,$$
where $T_n^+(t) := \max(T_n(t),0)$ and $T_n^-(t) := \max(-T_n(t),0)\,.$
Observe that (4.4) implies that 
$$0 \leq \int_{0}^{2\pi}{T_n^-(t) \, dt} = \int_{0}^{2\pi}{T_n^+(t) \, dt} \leq 2\pi \delta_n\,. \tag 4.5$$

Using the Bernstein inequality in $L_1$ for trigonometric polynomials $Q_n$ of degree $n$ defined by 
$$Q_n(t) := \delta_n - T_n(t) \geq 0\,, \qquad t \in {\Bbb R}\,, \tag 4.6$$
we have
$$\int_{0}^{2\pi}{|Q_n^{\prime}(t)| \, dt} \leq n \int_0^{2\pi}{|Q_n(t)| \, dt} = n \int_0^{2\pi}{Q_n(t) \, dt} \leq 2\pi n\delta_n\,.$$
Hence, recalling (4.1), (4.5), (4.6), we get
$$\int_{0}^{2\pi}{|T_n^{\prime}(t)| \, dt} \leq 2\pi n\delta_n\,. \tag 4.7$$
Let $k \geq 1$ be the largest integer such that $2^k \delta_n < (n+\delta_n)/2$. We define the sets
$$\eqalign{B_1 := & \{t \in [0,2\pi]: -\delta_n <  T_n(t) \leq \delta_n\}\,, \cr 
           A_j := & \{t \in [0,2\pi]: - 2^j\delta_n <  T_n(t) \leq - 2^{j-1}\delta_n\}\,, \qquad j=1,2,\ldots,k\,, \cr 
           B 2 := & \{t \in [0,2\pi]: -n \leq T_n(t) \leq -n/4\}\,. \cr}$$
Observe that 
$$[0,2\pi] =  B_1 \cup B_2 \cup_{j=1}^k{A_j}\,, \tag 4.8$$
$$m(B_1) \leq 2\pi\,, \tag 4.9$$
$$m(A_j)2^{j-1}\delta_n \leq \int_{A_j}{(T_n^-(t)) \, dt} \leq \int_0^{2\pi}{T_n^-(t) \, dt} \leq 2\pi \delta_n \,,$$ 
and hence
$$m(A_j) \leq 4\pi 2^{-j}\,, \qquad j=1,2,\ldots,k\,, \tag 4.10$$
where $m(E)$ denotes the linear Lebesgue measure of the set $E \subset {\Bbb R}$.
Observe also that 
$$m(B_2)n//4 \leq \int_{B_2}{T_n^-(t) \, dt} \leq \int_{0}^{2\pi}{T_n^-(t) \, dt} \leq 2\pi \delta_n \,,$$
and hence
$$m(B_2) \leq 8\pi \frac{\delta_n}{n}\,. \tag 4.11$$

Applying the Bernstein-Szeg\H o inequality (Lemma 3.2) with `$m := -n$, and $M := \delta_n$, and recalling (4.4), we get
$$\left|T_n^{\prime}(t)\right|^2 + n^2\left|T_n(t) - \frac{M+m}{2}\right|^2 \leq n^2\left(\frac{M-m}{2}\right)^2\,, \qquad t \in {\Bbb R}\,, \tag 4.12$$

If $t \in B_1$, then (4.11) implies
$$\left|T_n^{\prime}(t)\right|^2 + n^2\left| -\delta_n - \delta_n/2 + n/2 \right|^2 \leq n^2\left(n/2 + \delta_n/2\right)^2\,.$$
Hence for $t \in B_1$ we have
$$\split \left|T_n^{\prime}(t)\right|^2 \leq & n^2\left( \left(n/2 + \delta_n/2) \right)^2 - n^2\left( n/2 - 3\delta_n/2 \right)^2 \right) \cr
= & n^2 (n - \delta_n)2\delta_n \leq 2n^3\delta_n\,, \cr \endsplit$$
from which
$$\left|T_n^{\prime}(t)\right| \leq 2^{1/2}n^{3/2}\delta_n^{1/2}\,, \qquad t \in B_1\,, \tag 4.13$$
follows.

If $t \in A_j, j=1,2,\ldots,k$, then (4.11) implies
$$\left|T_n^{\prime}(t)\right|^2 + n^2\left(-2^j \delta_n + n/2 - \delta_n/2 \right)^2 \leq n^2\left(n/2 + \delta_n/2\right)^2\,.$$
Hence for $t \in A_j, j=1,2,\ldots,k$, we have
$$\split \left|T_n^{\prime}(t)\right|^2 \leq & n^2\left( \left(n/2 + \delta_n/2 \right)^2 - \left(-2^j \delta_n + n/2 - \delta_n/2 \right)^2 \right) \cr 
= & n^2 \left( \delta_n + 2^j\delta_n \right) \left( n - 2^j \delta_n + \right)\,, \cr \endsplit$$
from which
$$\left|T_n^{\prime}(t)\right| \leq 2^{1/2}n^{3/2}2^{j/2}\delta_n^{1/2}\,, \qquad t \in A_j\,, \tag 4.14$$
follows.  

If $t \in B_2$, then (4.11) implies
$$\left|T_n^{\prime}(t)\right|^2 \leq n^2\left((n/2 + \delta_n/2)\right)^2\,,$$
from which
$$\left|T_n^{\prime}(t)\right| \leq n^2\,, \qquad t \in B_2\,. \tag 4.15$$

Let $A := \cup_{j=1}^k{A_j}$. Using (4.2), (4.3), (4.7) - (4.15)  we have

$$\split & \pi \frac{n(n+1)(2n+1)}{6} \leq \int_{0}^{2\pi}{\left|T_n^{\prime}(t)\right|^2 \, dt} \cr
\leq & \int_{B_1}{\left|T_n^{\prime}(t)\right|^2 \, dt} + \sum_{j=1}^k{\int_{A_j}{\left|T_n^{\prime}(t)\right|^2 \, dt}} + \int_{B_2}{\left|T_n^{\prime}(t)\right|^2 \, dt} \cr 
\leq & \left( \int_{0}^{2\pi}{|T_n^{\prime}(t)| \, dt} \right)  
\left(m(B_1) \max_{t \in B_1}{|T_n^{\prime}(t)|} + \sum_{j=1}^k{m(A_j) \max_{t \in A_j}{\left|T_n^{\prime}(t)\right|}} + m(B_2) \max_{t \in B_2}{|T_n^{\prime}(t)|} \right) \cr  
\leq & 2\pi n\delta_n \left( 2\pi 2^{1/2}n^{3/2}\delta_n^{1/2} + \sum_{j=1}^k{4\pi 2^{-j} 2^{1/2}n^{3/2}2^{j/2}\delta_n^{1/2}} + 8\pi \frac{\delta_n}{n}n^2 \right)\,, \cr 
= & 2\pi n^{5/2}\delta_n^{3/2} \left( 2\pi 2^{1/2}\pi + 4\pi 2^{1/2}\sum_{j=1}^k{2^{-j/2}} + 8\pi \delta_n^{1/2}n^{-1/2} \right) \cr 
\leq & 2\pi n^{5/2}\delta_n^{3/2} 2\pi \left( 1.42 + 6.84 + 4 \right) \cr \endsplit$$
hence, choosing $\delta_n := cn^{1/3}$, we get 
$$\frac 16 < c^{3/2} \pi 24.52\,.$$
Thus $\displaystyle{\left( \frac{1}{\pi 147.12} \right)^{2/3}} < c$, which is impossible if $\displaystyle{c := \frac{1}{60} < \left( \frac{1}{\pi 147.12} \right)^{2/3}}$.
\qed \enddemo

\demo{Proof of Corollary 2.2}
Let $T_n \in {\Cal LC}_n$. Applying Theorem 2.1 to the trigonometric polynomials $Q_n \in {\Cal LC}_n$ defined by $Q_n(t) := -T_n(t)$, 
we get the result,
\qed \enddemo

\demo{Proof of Remark 2.3}
Let $T_n \in {\Cal LC}_n$. By Parseval's formula we have
$$\int_{0}^{2\pi}{\left|T_n(t)\right|^2 \,  dt} = 2\pi \sum_{k=1}^{n}{1/2} \geq \pi n\,. \tag 4.16$$
Now let $T_n^+(t) := \max(t_n(t),0)$ and $T_n^-(t) := \max(-t_n(t),0)\,.$
Observe that $T_n(t)^2 := T_n^+(t)^2 + T_n^{-}(t)^2$ and by (4.16) we have either 
$$\int_{0}^{2\pi}{T_n^+(t)^2 \, dt} \geq \pi n/2$$ 
or 
$$\int_{0}^{2\pi}{T_n^-(t)^2 \, dt} \geq \pi n/2$$ and the statement of the remark follows.
\qed \enddemo

\enddocument